\documentclass[11pt,reqno]{amsart}

\title{On some families of modules for the current algebra}
\author{Mathew Bennet, \; Rollo Jenkins}
\address{Imecc -
Unicamp, Departamento de Matem\'{a}tica. Rua S\'{e}rgio Buarque de Holanda,
651, Cidade Universit\'{a}ria Zeferino Vaz. 13083-859 Campinas - SP, Brasil.
E-mails: mbenn002@gmail.com, rollojenkins@googlemail.com .}
\thanks{The first author was supported by Fapesp grant {\color{red}2012/06923-0} and the second author was supported by Fapesp grant {\color{red}2013/17654-3}. }

\usepackage{latexsym}
\usepackage{multicol}
\usepackage{verbatim}
\usepackage{amsmath, amscd}
\usepackage{color}
\usepackage[latin1]{inputenc}
\usepackage{enumerate}
\usepackage{amsfonts,amsmath,amsthm,amssymb,amscd}
\usepackage{latexsym}
\usepackage{enumerate}
\usepackage{mathrsfs}
\usepackage{pxfonts}
\usepackage[table,dvipsnames]{xcolor}
\definecolor{rollo}{rgb}{0,0.528,0}
\usepackage[%
	bookmarks=true,			
	unicode=true,			
	pdftitle={},		%
	pdfauthor={},	%
	pdfkeywords={} {},	
	colorlinks=true,		
	linkcolor=Red,			
	citecolor=rollo,		
	filecolor=magenta,		
	urlcolor=RoyalBlue			
  ]{hyperref}				

\definecolor{rollo}{rgb}{0,0.528,0}

\advance\textwidth by 1.2in \advance\oddsidemargin by -.6in \advance\evensidemargin by -.6in
\newtheorem*{cor}{Corollary}
\newtheorem*{lem}{Lemma}

\newtheorem*{prop}{Proposition}

\theoremstyle{definition}

\theoremstyle{definition}
\newtheorem*{thm}{Theorem}

\newtheorem*{coro}{Corollary}

\newcommand{\bydef}{\mathrel{\mathop:}=}

\newenvironment{pf}{\proof}{\endproof}
\newcounter{cnt}
 \makeatletter
\def\mydggeometry{\makeatletter\dg@YGRID=1\dg@XGRID=20\unitlength=0.003pt\makeatother}
\makeatother \theoremstyle{remark}

\numberwithin{equation}{section}

\let\bwdg\bigwedge
\def\bigwedge{{\textstyle\bwdg}}

\begin{document}

\newcommand{\thmref}[1]{Theorem~\ref{#1}}
\newcommand{\secref}[1]{Section~\ref{#1}}
\newcommand{\lemref}[1]{Lemma~\ref{#1}}
\newcommand{\propref}[1]{Proposition~\ref{#1}}
\newcommand{\corref}[1]{Corollary~\ref{#1}}
\newcommand{\remref}[1]{Remark~\ref{#1}}
\newcommand{\defref}[1]{Definition~\ref{#1}}
\newcommand{\er}[1]{(\ref{#1})}
\newcommand{\id}{\operatorname{id}}
\newcommand{\ord}{\operatorname{\emph{ord}}}
\newcommand{\sgn}{\operatorname{sgn}}
\newcommand{\Sgn}{\operatorname{\mathbf{sgn}}}
\newcommand{\Triv}{\operatorname{\mathbf{triv}}}
\newcommand{\wt}{\operatorname{\mathrm{wt}}}
\newcommand{\tensor}{\otimes}
\newcommand{\from}{\leftarrow}
\newcommand{\nc}{\newcommand}
\newcommand{\rnc}{\renewcommand}
\newcommand{\dist}{\operatorname{dist}}
\newcommand{\qbinom}[2]{\genfrac[]{0pt}0{#1}{#2}}
\nc{\cal}{\mathcal} \nc{\goth}{\mathfrak} \rnc{\bold}{\mathbf}
\renewcommand{\frak}{\mathfrak}
\newcommand{\supp}{\operatorname{supp}}
\newcommand{\Irr}{\operatorname{Irr}}
\renewcommand{\Bbb}{\mathbb}
\nc\bomega{{\mbox{\boldmath $\omega$}}} \nc\bpsi{{\mbox{\boldmath $\Psi$}}}
 \nc\balpha{{\mbox{\boldmath $\alpha$}}}
 \nc\bpi{{\mbox{\boldmath $\pi$}}}
\nc\bsigma{{\mbox{\boldmath $\sigma$}}} \nc\bcN{{\mbox{\boldmath $\cal{N}$}}} \nc\bcm{{\mbox{\boldmath $\cal{M}$}}} \nc\bLambda{{\mbox{\boldmath
$\Lambda$}}}

\newcommand{\lie}[1]{\mathfrak{#1}}

\newcommand{\tlie}[1]{\tilde{\mathfrak{#1}}}
\newcommand{\hlie}[1]{\hat{\mathfrak{#1}}}
\newcommand{\tscr}[1]{\tilde{\mathscr{#1}}}
\newcommand{\hscr}[1]{\hat{\mathscr{#1}}}
\newcommand{\hcal}[1]{\hat{\mathcal{#1}}}
\newcommand{\tcal}[1]{\tilde{\mathcal{#1}}}

\makeatletter
\def\section{\def\@secnumfont{\mdseries}\@startsection{section}{1}%
  \z@{.7\linespacing\@plus\linespacing}{.5\linespacing}%
  {\normalfont\scshape\centering}}
\def\subsection{\def\@secnumfont{\bfseries}\@startsection{subsection}{2}%
  {\parindent}{.5\linespacing\@plus.7\linespacing}{-.5em}%
  {\normalfont\bfseries}}
\makeatother
\def\subl#1{\subsection{}\label{#1}}
 \nc{\Hom}{\operatorname{Hom}}
  \nc{\mode}{\operatorname{mod}}
\nc{\End}{\operatorname{End}} \nc{\wh}[1]{\widehat{#1}} \nc{\Ext}{\operatorname{Ext}} \nc{\ch}{\text{ch}} \nc{\ev}{\operatorname{ev}}
\nc{\Ob}{\operatorname{Ob}} \nc{\soc}{\operatorname{soc}} \nc{\rad}{\operatorname{rad}} \nc{\head}{\operatorname{head}}
\def\Im{\operatorname{Im}}
\def\gr{\operatorname{gr}}
\def\chg{\operatorname{ch_{gr}}}
\def\ch{\operatorname{ch}}
\def\mult{\operatorname{mult}}
\def\Max{\operatorname{Max}}
\def\ann{\operatorname{Ann}}
\def\sym{\operatorname{sym}}
\def\loc{\operatorname{loc}}
\def\coin{\operatorname{coin}}
\def\dec{\operatorname{Dec}}
\def\maj{\operatorname{Maj}}
\def\stab{\operatorname{STab}} 
\def\tab{\operatorname{Tab}}                    
\def\part{\operatorname{\mathcal P}}
\def\Res{\operatorname{\br^\lambda_A}}
\def\und{\underline}
\def\Lietg{$A_k(\lie{g})(\bsigma,r)$}

 \nc{\Cal}{\cal} \nc{\Xp}[1]{X^+(#1)} \nc{\Xm}[1]{X^-(#1)}
\nc{\on}{\operatorname} \nc{\Z}{{\bold Z}} \nc{\J}{{\cal J}} \nc{\CC}{{\bold C}} \nc{\Q}{{\bold Q}}
\renewcommand{\P}{{\cal P}}
\nc{\N}{{\Bbb N}} \nc\boa{\bold a} \nc\bob{\bold b} \nc\boc{\bold c} \nc\bod{\bold d} \nc\boe{\bold e} \nc\bof{\bold f} \nc\bog{\bold g}
\nc\boh{\bold h} \nc\boi{\bold i} \nc\boj{\bold j} \nc\bok{\bold k} \nc\bol{\bold l} \nc\bom{\bold m} \nc\bon{\bold n} \nc\boo{\bold o}
\nc\bop{\bold p} \nc\boq{\bold q} \nc\bor{\bold r} \nc\bos{\bold s} \nc\boT{\bold t} \nc\boF{\bold F} \nc\bou{\bold u} \nc\bov{\bold v}
\nc\bow{\bold w} \nc\boz{\bold z} \nc\boy{\bold y} \nc\ba{\bold A} \nc\bb{\bold B} \nc\bc{\bold C} \nc\bd{\bold D} \nc\be{\bold E} \nc\bg{\bold
G} \nc\bh{\bold H} \nc\bi{\bold I} \nc\bj{\bold J} \nc\bk{\bold K} \nc\bl{\bold L} \nc\bm{\bold M} \nc\bn{\bold N} \nc\bo{\bold O} \nc\bp{\bold
P} \nc\bq{\bold Q} \nc\br{\bold R} \nc\bs{\bold S} \nc\bt{\bold T} \nc\bu{\bold U} \nc\bv{\bold V} \nc\bw{\bold W} \nc\bz{\bold Z} \nc\bx{\bold
x} \nc\KR{\bold{KR}} \nc\rk{\bold{rk}} \nc\het{\text{ht }}

\nc\toa{\tilde a} \nc\tob{\tilde b} \nc\toc{\tilde c} \nc\tod{\tilde d} \nc\toe{\tilde e} \nc\tof{\tilde f} \nc\tog{\tilde g} \nc\toh{\tilde h}
\nc\toi{\tilde i} \nc\toj{\tilde j} \nc\tok{\tilde k} \nc\tol{\tilde l} \nc\tom{\tilde m} \nc\ton{\tilde n} \nc\too{\tilde o} \nc\toq{\tilde q}
\nc\tor{\tilde r} \nc\tos{\tilde s} \nc\toT{\tilde t} \nc\tou{\tilde u} \nc\tov{\tilde v} \nc\tow{\tilde w} \nc\toz{\tilde z}

\begin{abstract} Given a finite-dimensional module, $V$, for a finite-dimensional, complex, semi-simple Lie algebra $\lie g$ and a positive integer $m$, we construct a family of graded modules for the current algebra $\lie g[t]$ indexed by simple $\CC\lie S_m$-modules.  These modules have the additional structure of being free modules of finite rank for the ring of symmetric polynomials and so can be localized to give finite-dimensional graded $\lie g[t]$-modules. We determine the graded characters of these modules and show that if $\lie g$ is of type $A$ and $V$ the natural representation, these graded characters admit a curious duality.  
\end{abstract}

\maketitle

\section{Introduction}

The current algebra associated to a simple Lie algebra $\lie g$ is the Lie algebra $\lie g[t] = \lie g\otimes \bc[t]$ with the obvious bracket (see Section \ref{def:current algebra}).  
The study of graded representations for the current algebra $\lie g[t]$ has been of significant interest for several decades.  One of the reasons for this is that the current algebra is a maximal parabolic subalgebra of the associated affine Kac-Moody Lie algebra and many interesting representations for the affine algebra specialize to graded modules for the current algebra.  One can also get modules by taking the so called graded-limit of modules for the quantum affine algebras (see for example \cite{C} and \cite{M}).

Of particular interest is the category $\cal I$, comprising graded modules $M$ for the current algebra with the condition that the graded components are finite dimensional. The simple modules in $\cal I$ are indexed by pairs, $\{(\lambda, r)\in P^+\times \bz\}$, where $P^+$ is the set of dominant integral weights for the Lie algebra $\lie g$.  The category $\cal I$ is not a semi-simple category; there are many indecomposable, yet reducible, modules in $\cal I$.  The subcategory of $\cal I$ consisting of objects with only finitely many non-zero graded components of negative grade is known to admit a BGG-style reciprocity (see \cite{CI} for the most general case) and a tilting theory (see \cite{BC}), while the subcategory consisting of finitely-generated modules is a motivating example of an affine highest weight category (see \cite{Kh} and \cite{Kl}).

Our work can be thought of as a generalization of the following classical construction. Because we will work exclusively over $\CC$, the representation theory of the symmetric group $\lie S_m$ is semi-simple, with the simple modules indexed by partitions, $\gamma\in\part(m)$, of $m$. Given a finite-dimensional $\lie g $-module $V$, the module $N= V^{\otimes m}$ is naturally a bimodule for $\lie g$ and $\lie S_m$ (both with left actions) and the actions commute. If we decompose $N \cong_{\lie S_m} \bigoplus_\gamma S(\gamma) \otimes \Hom_{\lie S_m}(S(\gamma), N)$ as a $\CC\lie S_m$-module, we have an action of $\lie g$ on the multiplicity spaces $\Hom_{\lie S_m}(S(\gamma),N)$.  As special cases, the multiplicity space associated to the trivial module is the space of invariants $N^{\lie S_m}$, while that associated to the sign representation is the exterior power $\wedge^m V$.

In our situation, we start with a finite-dimensional $\lie g$-module $V$, and let $V\otimes \bc[t]$ be a $\lie g[t]$-module with the natural action.  The tensor product $M = (V\otimes \bc[t])^{\otimes m}$ is naturally a $\lie g[t]$-$\CC\lie S_m$-bimodule and again, the actions commute. The grading on $\bc[t]$ by powers of $t$ induces a grading on $M$.  Clearly, $\lie S_m$ preserves the finite-dimensional graded components of $M$ and so we can decompose $M$ as follows.
\[ M\cong_{\lie S_m}  \bigoplus_{k\in \bz}  \bigoplus_{\gamma\in \part(m)} S(\gamma) \otimes \Hom_{\lie S_m}(S(\gamma), M[k]).\]
These multiplicity spaces are naturally $\lie g[t]$-modules; we denote them
\[ B(\gamma, V)\bydef \bigoplus_{k\in \bz} \Hom_{\lie S_m}(S(\gamma), M[k]) .\]
Note that these modules lie in $\cal I$. The purpose of this paper is to study the modules $B(\gamma, V)$. 

We note here that the modules $B(\gamma, V)$ have additional structure.  Using the natural isomorphism $M \cong V^{\otimes m} \otimes \bc[t_1, \ldots, t_m]$, we see that $M$ admits a right action by the polynomial ring $\bc[t_1, \ldots, t_m]$, which is simply right multiplication.  This action commutes with the action of $\lie g[t]$, but not with the action of the symmetric group.  However, if we restrict this right action to the ring of symmetric polynomials $\ba_m$, then all three actions commute.  In fact, $M$ is a free $\ba_m$-module of finite rank. It follows from the celebrated result of Quillen and Suslin that $B(\gamma, V)$ is a free right $\ba_m$-module.  We let $\bold I_m$ be the unique graded maximal ideal in $\ba_m$ and define the graded localization $B_{\loc}(\gamma, V) := B(\gamma, V) \otimes _{\ba_m} \ba_m / \bi_m$.  This is a finite-dimensional graded $\lie g[t]$-module.  The main result of our paper is the following theorem.

\begin{thm} The multiplicity space $B(\gamma, V)$ is a graded module for the current algebra and so is $B_{\loc}(\gamma, V)$. The graded characters of these modules behave in the following way.
\begin{enumerate}[(i)]
\item For each $\gamma\in \part(m)$, the graded character of $B_{\loc}(\gamma, V)$ is given by 
\[
\chg B_{\loc}(\gamma, V) = \sum_{\mu \in P^+}\sum_{\sigma,\tau\in \part(m)} s_\mu(\tau, V) c_{\tau, \sigma}^\gamma f_\sigma(u) e(\cal O(\mu)).
\]
The notation here is defined in Section \ref{Section:Notation}.
\item Up to a fixed graded shift, the graded characters of localizations corresponding to conjugate partitions and dual $\lie g$-modules are dual to one another in the following way
\[
\chg B_{\loc}(\gamma, V) = u^{\binom{m}{2}} \chg B_{\loc}(\gamma^\vee, V^*)^*.
\]
\end{enumerate}
\end{thm}

Part (i) can easily be used to give the graded character of $\chg B(\gamma, V)$. When $V=V(\omega_1)$ is the natural representation, we will show that the coefficients $s_\mu(\tau, V)$ appearing in the formula above are Kostka numbers.

Some of these modules have already been studied in other contexts.  When $\lie g =\lie{sl}_{n+1}$ and  $V=V(\omega_1)$, the $\lie g[t]$-module $V\otimes \bc[t]$ is the global Weyl module of highest weight $\omega_1$.  The global Weyl module of highest weight $m\omega_1$ can be realized as the $\lie S_m$ invariants in $W(\omega_1)^{\otimes m}$ (\cite{CPWeyl} for $\lie {sl}_2$ and  \cite{FL} for $\lie{sl}_{n+1}$); that is, $B(\Triv, V(\omega_1))$.

In the case that $\lie g= \lie {sl}_2$, it is shown in \cite{BC2} that $W(1)^{\otimes m}$ admits a filtration by global Weyl modules (as does any summand) and that the tilting module of highest weight $m$ is (up to a grade shift) the exterior product of $W(\omega_1)$; that is $B(\Sgn, V(1))$. It follows from the first result that 
$\{ B(\tau, V(1)) : \tau \in \part(m), \; m \in \bz_{>0} \}$ is a two parameter family of modules admitting filtrations by global Weyl modules. The main theorem also gives a new formula for the graded character of the global Weyl modules $W(m \omega_1)$ and the tilting modules $T(m)$. We also note that the character equality implied by \cite{BC2} Proposition 3.7 can be thought of as a special case of part (ii) of the Theorem.

We also note the following generalizations of this construction.  Given a commutative associative algebra A , $\lie g \otimes A$ is naturally a Lie algebra.  Given a $\lie g$-module $V$, the space $V\otimes A$ is naturally a $\lie g\otimes A$-module, and $(V\otimes A)^{\otimes m}$ is a bimodule for $\lie g\otimes A$ and $\lie S_m$ such that the actions commute.  Therefore the multiplicity spaces associated to simple modules $S(\tau)$ will be $\lie g\otimes A$-modules.  A key component to understanding the multiplicity spaces when $A = \bc[t]$ is understanding the $\lie S_m$ structure of $\bc[t]^{\otimes m} \cong \bc[t_1, \ldots, t_m]$.  However, very little can be said about the $\lie S_m$ structure of $A^{\otimes m}$ in this generality. 

\textit{Acknowledgements: The authors would like to thank Adriano de Moura for many useful discussions.}

\section{Notation}\label{Section:Notation} Throughout this paper we will work over the complex numbers. All unmarked tensor products are understood to be over $\CC$. All graded algebras will be positively graded over the integers. 

\subsection{Lie theory, characters and gradings}
Let $\lie g$ denote a finite-dimensional complex simple Lie algebra, $\lie h \subset \lie g$ a fixed Cartan subalgebra of rank $n$.  We set $I = \{1, \ldots, n\}$ and let $\{ \alpha_i \, : \, i\in I\}$ be a set of simple roots for $\lie h^*$ with respect to $\lie h$.  We let $R$ denote the set of roots, $R^+$ the positive roots, $Q$ the root lattice and $Q^+$ the positive root lattice.  We also let $P$ be the set of weights and $P^+$ the dominant integral weights, with $\{ \omega_i\; : \; i\in I \}$ the fundamental weights.  We put a partial order on $\lie h^*$ by letting $\lambda \ge \mu$ if $\lambda - \mu \in Q^+$.  

Let $\cal W \subset \mathrm{Aut}(\lie h^*)$ be the Weyl group and $w_0$ the longest word.  For $\alpha \in R$ let $\lie g_\alpha$ be the corresponding root space.  We have a decomposition $\lie g = \left(\oplus_{\alpha \in R} \lie g_\alpha\right)\oplus\lie h $ and set $\lie n^{\pm} = \oplus_{\alpha \in R^+} \lie g_{\pm \alpha}$.  Let $\{ x_\alpha,\; h_i \, : \, \alpha \in R^+,\;i \in I \}$ be a Chevalley basis for $\lie g$. We set $x^{\pm}_i = x^{\pm}_{\alpha_i}$.

Given a $\lie g$-module $M$ and $\lambda \in \lie h^*$, denote by $M_\lambda$ the $\lambda$ weight space of $M$; that is $M_\lambda \bydef \{ m \in M \, :\,  h\cdot m = \lambda(h)m\text{, for all } h \in \lie h \}$. Finite-dimensional $\mathfrak g$-modules decompose as a direct sum of their weight spaces. Any module with such a decomposition is referred to as a \textbf{weight module}. The finite-dimensional simple $\lie g$-modules each correspond to a dominant integral weight $\lambda \in P^+$, and we call this module $V(\lambda)$; it is generated by a vector $v_\lambda$ of weight $\lambda$ such that $\lie n^+\cdot v_\lambda =0$.  The finite dimensional representation theory of $\lie g$ is semi-simple, meaning that every finite dimensional $\lie g$-module is isomorphic to a direct sum of simple modules.  

Given a finite dimensional $\lie g$-module $V$, the dual space $V^*$ is also a finite dimensional $\lie g$-module, and hence a weight module. Note that 
\begin{equation} (V_\lambda)^* \cong  (V^*)_{-\lambda}. \end{equation}
It is well known that if $\lambda\in P^+$ then $\lambda^\vee := -\omega_0\cdot \lambda \in P^+$. On simple modules, this duality satisfies $V(\lambda)^*\cong V(\lambda^\vee)$.

The \textbf{character} of a finite-dimensional $\lie g$-module, $M$, is the formal sum in the group ring of the weight lattice,
\[ \ch_{\lie g} M \bydef \sum_{\lambda \in P} \dim M_\lambda e(\lambda) \in \bz[P].\]
A standard result says that if $w\in \cal W$ and $\lambda \in \lie h^*$, then $\dim M_\lambda = \dim M_{w\cdot\lambda}$ for any finite-dimensional module $M$.  
If we let $\cal O(\lambda)= \{ w\cdot\lambda : w\in \cal W \}$ and $e(\cal O(\lambda)) = \sum_{\mu \in \cal O(\lambda)} e(\mu)$ then we can also write the character of $M$ as
\[ \ch_{\lie g} M = \sum_{\lambda \in P^+} \dim M_\lambda e(\cal O(\lambda)).\]

Given a finitely generated, graded algebra $A=\oplus_{i\in \bz}A[i]$, we call an $A$-module $M$ a \textbf{graded} module if it decomposes as a vector space into graded pieces $M = \oplus_{k\in \bz} M[k]$ such that $A[i] M[j] \subset M[i+j]$.  If each graded component is finite dimensional we can associate to $M$ its \textbf{Hilbert series},
\[
H(M) \bydef \sum_{k\in \bz} \dim M[k] u^k \in \bz[[ u, u^{-1}]]
\] where $u$ is an indeterminate.

If $M$ is also a $\lie g$-module, we define its \textbf{graded character} to be the following power series with coefficients in the group ring $\bz[P]$;
\[ \chg M \bydef \sum_{r\in \bz} \ch_{\lie g} M[r] u^r \in \bz[P][[u^{\pm 1}]].\]

\subsection{The coinvariant ring for $\mathfrak S_m$}  

Let $\part(m)= \{ \tau= (\tau_1 \ge \tau_2 \ge \cdots \ge \tau_m) \, :\, \sum \tau_i =m \}$ be the set of partitions of $m$.  Given $\tau \in \part(m)$ we let $\tau^\vee$ denote the conjugate partition. The complex irreducible representations of $\lie S_m$ are indexed by partitions of $m$. Given such a partition, $\tau$, we will denote the corresponding simple $\CC\lie S_m$-module $S(\tau)$.  We will denote the one-dimensional modules corresponding to the trivial and sign representations by $\Triv$ and $\Sgn$ respectively.  The group algebra $\bc \lie S_m$ has a comultiplication induced by the assignment $g \mapsto g \otimes g$;  in this way, we define the tensor product of $\CC\lie S_m$-modules.  We can also define the dual representation $S^* = \Hom_{\lie S_m}(S, \bc)$, which on simple modules satisfies $S(\gamma)^* \cong S(\gamma)$. 

Define non-negative integers, $c_{\tau, \sigma}^\gamma$, by 
\begin{equation}\label{coeff1} S(\tau)\otimes S(\sigma) \cong \bigoplus_{\gamma\in \part(m)} c_{\tau, \sigma}^\gamma S(\gamma); \end{equation}
the $c_{\tau, \sigma}^\gamma$ are known as Kronecker coefficients. Note that the decomposition above depends on a choice of basis for each isotypic component of $S(\tau)\otimes S(\sigma)$ and hence is not canonical.  
However, for any finite-dimensional module, $N$, the following decomposition \emph{is} canonical.
\begin{equation}\label{candecomp} N \cong \bigoplus_{\tau\in \part(m)} S(\tau) \otimes \Hom_{\lie S_m}(S(\tau), N) \end{equation}
where the dimension of $\Hom_{\lie S_m}(S(\tau), N)$ records the multiplicity of the isotypic component corresponding to $\tau$. 

Given a finite dimensional $\lie g$ module $V$, the module $N = V^{\otimes m}$ is a left $\lie S_m$-module via permutation of tensorands, and this action commutes with the left action of $\lie g$.  In particular $\lie S_m$ preserves the weight spaces of $N$.  By the discussion above, for each $\mu \in P$, the weight space $N_\mu $ decomposes into a direct sum of simple $\lie S_m$-modules.  We define non-negative integers $s_\mu(\tau, V)$ to satisfy
\begin{equation}\label{coeff2} N_\mu \cong_{\lie S_m} \bigoplus_\gamma s_\mu(\gamma, V) S(\gamma). \end{equation}

Let $A_m \bydef \bc[t_1, \ldots, t_m]$, the polynomial ring in $m$ indeterminates, which we consider to be $\bz_{\ge 0}$ graded in the natural way.  Then $\lie S_m$ acts on $A_m$ by permuting the $t_i$.  Clearly, $\lie S_m$ preserves the graded components of $A_m$, so $A_m$ is a graded $\CC\lie S_m$-module.  Let $\bold A_m \bydef A_m^{\lie S_m}$ be the ring of polynomials invariant under the action of $\lie S_m$ and $(\bold A_m)_+ = \oplus_{k\ge 1} \bold A_m[k]$.  It is know that $\ba_m$ is itself a polynomial ring, and that $A_m$ is a free module over $\ba_m$.  The \textbf{coinvariant ring}\footnote{also known as the co(in)variant algebra} is the quotient 
\[A_m^{\coin} \bydef A_m / (\bold A_m)_+.\]  

It is well known (for example, \cite{Chevalley}) that the module $A_m^{\coin}$ is isomorphic to the regular representation of $\lie S_m$; however, in the category of graded representations, these two are not isomorphic: the coinvariant has a non-trivial grading which reflects some of the more subtle features of the symmetric group, such as the ordering of simple $\lie S_m$-modules induced by the natural ordering of partitions. More precisely, given a partition $\gamma\in\part(m)$, we can decompose the $S(\gamma)$-isotypic component of $A_m^{\coin}$ into graded pieces. We denote the Hilbert series of these components by $f_\gamma(u)$; that is to say,
\begin{equation}\label{coeff3} H(A_m^{\coin}) = \sum_{\gamma \in \part(m)} f_\gamma(u) \dim S(\gamma). \end{equation}
Note that the $f_\gamma(u)$ are actually polynomials.

The proof of following result on the Hilbert series of free $\ba_m$ modules is straight forward
\begin{lem}\label{freemod} Let $M$ be a free, finitely generated, graded $\ba_m$ module with graded basis $\cal S$ such that there are $n_r$ elements of $\cal S$ of grade $r$. Then $H(M) = H(\ba_m) \sum n_r u^r$.
\end{lem}

\subsection{The Current Algebra}\label{def:current algebra}
For a Lie algebra, $\lie a$, the \textbf{current algebra of $\lie a$} is the vector space $\lie a[t] = \lie a\otimes \bc[t]$ with the Lie bracket determined by the rule $[a\otimes t^r,a'\otimes t^s]\bydef [a,a']\otimes t^{r+s}$ for $a,a'\in \lie a$. Let $\bu(\lie a)$ denote the universal enveloping algebra of $\lie a$.  The current algebra is graded by putting $t$ in degree one, as is $\bu(\lie a[t])$.  Note that $\lie g\otimes 1 \subset \lie g[t]$ is isomorphic to $\lie g$ and we will simply write $\lie g \subset \lie g[t]$. The universal enveloping algebra $\bu(\lie a)$ is an associative Hopf algebra with co-multiplication induced by the assignment $\Delta(x) = 1\otimes x + x \otimes 1$ for $x\in\lie a$.  This is a map of graded algebras for the case $\bu(\lie a[t])$.

We let $\cal I$ be the category of graded $\lie g[t]$-modules, $M$, such that each graded component satisfies $\dim M[k] < \infty$ with graded morphisms (that is, they preserve degrees). Note that, because the graded components are finite-dimensional $\lie g$-modules, the objects in $\cal I$ are necessarily weight modules. The simple modules in $\cal I$ are indexed by $(\lambda, r)\in P^+ \times \bz$; the simple module $V(\lambda, r)$ is isomorphic as a $\lie g$-module to $V(\lambda)$, and is concentrated in the $r^{th}$ graded component.

The category $\cal I$ is not strictly a tensor category:  given two objects $M$ and $N$, the tensor product $M \otimes N$ satisfies $(M\otimes N)[k] = \sum_{i\in \bz} M[k-i]\otimes N[i]$.  This will be an object of $\cal I$ if and only if this sum is finite for all $k$.  A sufficient condition is that both $M$ and $N$ have a lower bound on their grades. 

Given an object, $N\in\cal I$, its \textbf{graded dual} is the module $N^*=\oplus_{k\in \bz} \Hom(N[k], \bc)$, which has graded components
\[
(N^*)[k] = (N[-k])^*.
\] 
The following lemma explains the connection between the graded characters of a module $N\in \cal I$ and that of its dual, $N^*$.

\begin{lem}\label{duallem} Let $\chg N = \sum_{\lambda\in P^+} g_\lambda(u) e(\cal O(\lambda))$, where $g_\lambda(u)\in \bz[[u^{\pm 1}]]$.  Then the graded character of $N^*$ is given by 
\[ \chg N^* = \sum_{\lambda \in P^+} g_\lambda(u^{-1}) e(\cal O(\lambda^\vee)).\]
\end{lem}
\begin{pf} It is enough to show that the coefficient of $e(\cal O(\lambda^\vee)$ is $g_\lambda(u^{-1})$, for which it is enough to show that the coefficient of $(N^*)_{\lambda^\vee}$ is $g_\lambda(u^{-1})$. For this, it is enough to show that the dimension of $(N^*)_{\lambda^\vee}[-k]$ is equal to the dimension of $N_\lambda[k]$.

If $N_\lambda \ne 0$ then $(N_\lambda)^* = (N^*)_{-\lambda} \ne 0$.  By the invariance of dimensions of weight spaces under the action of the Weyl group and equation \ref{dimdual} it follows that 
\[  \dim N[k]_\lambda = \dim N[k]_{\omega_0 \lambda} = \dim (N^*)[-k]_{\lambda^\vee }.\]   
The result follows.
\end{pf}

The subcategory of $\cal I$ where objects have a lowest graded component admits a tilting theory(\cite{BC}); for each $(\lambda, r)\in P^+\times \bz$ there is an indecomposable tilting module $T(\lambda, r)$ which admits a filtration by standard modules and a filtration by costandard modules. The standard modules are the global Weyl modules $W(\lambda, r)$: these are universal highest weight modules, and can be defined using generators and relations.  If $\lambda \ne 0$ the global Weyl module is infinite dimensional. The costandard objects are the graded duals of so called local Weyl modules $(W_{\loc}(\lambda, r))^*$. The local Weyl module is a finite dimensional quotient of the global Weyl module, and can be defined using generators and relations.  Global and local Weyl modules were introduced in \cite{CPWeyl} and defined in broad generality in \cite{CFK}.  The modules $W(\lambda, r)$ and $W(\lambda, s)$ are grade shifts of each other.

If $V$ is a finite-dimensional $\lie g$-module then $V[t] \bydef V\otimes \bc[t]$ is an object in $\cal I$, with an action defined by letting $x\otimes t^r\cdot v\otimes t^s \bydef (x\cdot v) \otimes t^{r+s}$, for $x\in \lie g$, $v\in V$ and $r,s\in \bz$.  We have 
$\chg V[t] = \sum_{r\ge 0} \ch_{\lie g}V u^r$. 
In \cite{BCGM} it is shown that the fundamental global Weyl module $W(\omega_i, 0) \cong V(\omega_i)\otimes \bc[t]$.

\section{A natural construction}
For $V$ a finite-dimensional $\lie g$-module, the module $M = (V[t])^{\otimes m} \in \cal I$ admits a natural right action of $\lie S_m$, permuting tensorands; therefore, $M$ is a $\lie g[t]$-$\CC\lie S_m$-bimodule and the two actions commute. Because $\lie S_m$ preserves the finite-dimensional graded components, the decomposition of $M$ as a representation of $\lie S_m$ follows Equation \ref{candecomp}:
\begin{align*}
M  &\cong_{\lie S_m} \bigoplus_{\gamma\in \part(m)} \bigoplus_{k\in \bz} S(\gamma) \otimes \Hom_{\lie S_m}(S(\gamma), M[k]).
\end{align*}
For each $k\in \Z$, we set $B(\gamma, V)[k] = \Hom_{\lie S_m}(S(\gamma), M[k])$ and we define
\begin{align*}
B(\gamma, V) &= \bigoplus_{k\in \Z} \Hom_{\lie S_m}(S(\gamma), M[k]).
\end{align*}
We will show at the end of the section that $B(\gamma, V)$ is a $\lie g[t]$-module.  

Consider the vector space isomorphism 
\begin{equation}\label{mis}
M=(V\otimes \bc[t])^{\otimes m} \cong V^{\otimes m}\otimes A_m. 
\end{equation}

Using equation \ref{mis} we see that $M$ admits the structure of a right $A_m$-module.  This action does not commute with the action of $\lie S_m$; however, if we restrict this right action to $\ba_m$, the actions commutes with $\lie S_m$ and $\lie g[t]$.
Because the polynomial ring $A_m$ is a free graded $\ba_m$-module of rank $=\dim A_m^{\coin}$, we see that $M$ is a free $\bold A_m$-module of rank $= (\dim V)^m \times \dim A_m^{\coin}$. The $B(\gamma, V)$ are then projective $\bold A_m$ modules, and the Quillen-Suslin Theorem tells us that the $B(\gamma, V)$ are in fact free (graded) $\bold A_m$-module.  Given $\bold J$ a maximal ideal of $\ba_m$, we define the localization of $B(\gamma, V)$ at $\bold J$ to be
\[ B_{\bold J}(\gamma, V):= B(\gamma, V)\otimes_{\ba_m} \ba_m / \bold J\]
and analogously define $M_{\bold J}$. These will always be a finite dimensional vector space, but will only be a graded vector space if $\bold J$ is a graded ideal.  We let $\bold I_m$ be the unique graded maximal ideal in $\ba_m$, and denote $B_{\bold I_m}(\gamma, V)$ by $B_{\loc}(\gamma, V)$ and $M_{\bi_m}$ by $M_{\loc}$.

The main result of this paper is the following theorem.  See equations \ref{coeff1}, \ref{coeff2} and \ref{coeff3} for the definitions of the coefficients.

\begin{thm}\label{main} The multiplicity space $B(\gamma, V)$ is a graded module for the current algebra and so is $B_{\loc}(\gamma, V)$.  The graded characters of these modules behave in the following way. 
\begin{enumerate}[(i)]
\item For each $\gamma\in \part(m)$, the graded character of $B(\gamma, V)$ is given by 
\[
\chg B_{\loc}(\gamma, V) = \sum_{\mu \in P^+}\sum_{\sigma,\tau\in \part(m)} s_\mu(\tau, V) c_{\tau, \sigma}^\gamma f_\sigma(u) e(\cal O(\mu)).
\]
\item Up to a fixed graded shift, the graded characters of localizations corresponding to conjugate partitions and dual $\lie g$-modules are dual to one another as $\lie g[t]$-modules in the following sense 
\[
\chg B_{\loc}(\gamma, V) = u^{\binom{m}{2}} \chg B(\gamma^\vee, V^*)^*
\]
\end{enumerate}
\end{thm}

This theorem is proved in the next section. We end this section by stating the following, which is an immediate corollary of Theorem \ref{main}(i) and Lemma \ref{freemod} 
\begin{cor} The graded character of $B(\gamma, V)$ is given by
\[
\chg B(\gamma, V) =  H(\ba_m)\sum_{\mu \in P^+}\sum_{\sigma,\tau\in \part(m)} s_\mu(\tau, V) c_{\tau, \sigma}^\gamma f_\sigma(u) e(\cal O(\mu)).
\]
\end{cor}

We now explain why the (graded) multiplicity space $B(\gamma, V) = \Hom_{\lie S_m}(S(\gamma), M)$ is a $\lie g[t]$ module. Recall that the actions of $\lie g[t]$ and $\lie S_m$ on $M$ commute. This means that, as elements of $\End(M)$, we can view $\lie g[t] \subset \End_{\lie S_m}(M)$. It is now sufficient to show that $B(\gamma, V)$ is a module for $\End_{\lie S_m}(M)$. To the end, let $f\in B(\gamma, V)$ and $\phi \in \End_{\lie S_m}(M)$. Let $\phi\cdot f$ denote the composition $f$ followed by $\phi$. Then 
\[ S(\gamma) \stackrel{f}\rightarrow M \stackrel{\phi}\rightarrow M\] 
is clearly an element of $B(\gamma, V)$.

\section{Proof of Theorem \ref{main}}
In this section we will look more closely at several $\CC\lie S_m$-modules.  

\subsection{Hilbert polynomials for the coinvariant ring}
Recall that we defined the polynomials $f_\sigma(u)$ by the decomposition $H(A_m^{\coin})= \sum_\sigma f_\sigma(u) S(\sigma)$. In this section we prove the following proposition
\begin{prop}\label{thef} Let $\sigma\in \part(m)$ and $\sigma^\vee$ denote its conjugate. The polynomial $f_\sigma(u)$ satisfies
\[ f_\sigma(u) = u^{\binom{m}{2}} f_{\sigma^\vee}(u^{-1}).\] 
\end{prop}

The first step is to understand the polynomials $f_\sigma(u)$ better.  Let $\tab(\sigma)$  (respectively $\stab(\sigma)$) denote the set of tableau (respectively the set of standard tableau) of shape $\sigma$.   Given $T\in \tab(\sigma)$  we define its \textbf{descent} by $\dec(T) \bydef \{ a\, :\, a+1\text{ is in a row strictly below the row of }a\}$.  Define the \textbf{major index} of $T$ to be the integer $\maj_T = \sum_{a \in \dec(T)} a$.  For each $\sigma$ we define a polynomial 
\[ \maj_\sigma(u) = \sum_{T \in \stab(\sigma)} u^{\maj_T}.\]
The graded multiplicity of $S(\sigma)$ in $A_m^{\coin}$ is given by the major index, 
\begin{equation} f_\sigma(u) = \maj_\sigma(u); \end{equation} see, for example, \cite{Humphreys}. 

\begin{lem} The major indexes of the tableau $T$ and $T^\vee$ are related by $\maj_T + \maj_{T^\vee} = \binom{m}{2}$.  \end{lem}
\begin{pf} It is easy to see that the conjugate map induces a bijection of sets $\vee:\stab(\sigma) \to \stab(\sigma^\vee)$ sending $T\mapsto T^\vee$.

If $b\in\{1, \ldots, m-1\}$ is such that $b\notin \dec(T)$, then $b+1$ must be in the same row as $b$, and to the right. Then $b+1$ must be below $b$ in $T^\vee$, and so $b\in \dec(T^\vee)$.
It follows that $\dec(T) \cap \dec(T^\vee) = \emptyset$ and that $\dec(T) \cup \dec(T^\vee)= \{1, \ldots, m-1\}$.  Hence $\maj_T + \maj_{T^\vee} = \sum_{i=1}^{m-1} i = \binom{m}{2}$.
\end{pf}
Proposition \ref{thef} is an immediate corollary.

\subsection{The simple modules}
Recall the Kronecker coefficients are non-negative integers defined by $S(\tau)\otimes S(\sigma) = \oplus c_{\tau, \sigma}^\gamma S(\gamma)$.  We will show that
\begin{lem}\label{thec} The Kronecker coefficients satisfy $c_{\tau, \sigma}^\gamma = c_{\tau, \sigma^\vee}^{\gamma^\vee}$. \end{lem}
The key tool to proving this is the fact \cite[Exercise 4.51]{FultonHarris} that $S(\tau)\otimes \Sgn = S(\tau^\vee)$.

\begin{pf} We start with $S(\tau) \otimes S(\sigma) = \oplus c_{\tau, \sigma}^\gamma S(\gamma)$ and then tensor by $\Sgn$ on the right. This shows us that
\[ S(\tau) \otimes S(\sigma^\vee) = \bigoplus_\gamma c_{\tau, \sigma}^\gamma S(\gamma^\vee).\]
However, by definition we also have 
\[ S(\tau) \otimes S(\sigma^\vee) = \bigoplus_\gamma c_{\tau, \sigma^\vee}^{\gamma^\vee} S(\gamma^\vee).\]
Thus we must have $c_{\tau, \sigma}^\gamma = c_{\tau, \sigma^\vee}^{\gamma^\vee}$.
\end{pf}

\subsection{Weight modules}
\begin{prop}\label{wtspc} Let $V$ be a $\lie g[t] - \lie S_m$ bimodule and  $\mu \in \lie h^*$.  Then we have an isomorphism of $\lie S_m$-modules $V_\mu \cong_{\lie S_m} (V^*)_{-\mu}$. \end{prop}
\begin{pf} 
This follows from the fact that the $\lie g[t]$ duality satisfies $(V_\mu)^* \cong (V^*)_{-\mu}$ together with the fact that the $\lie S_m$ dual of a simple module is itself. 
\end{pf}

\subsection{Proof of the main result}
We can now prove Theorem \ref{main}.  We have already established the opening remark.  We will start by proving the first equality.  We freely use the notation of the previous sections.  Recall that for all $\mu \in P^+$, the weight space $(M_{\loc})_\mu$ is a $\lie S_m$-module, and 
\[ (M_{\loc})_\mu = (V^{\otimes m} \otimes A_m^{\coin})_\mu =  (V^{\otimes m})_\mu \otimes A_m^{\coin}. \]
Using the decompositions of $(V^{\otimes m})_\mu$ and $A_m^{\coin}$ as $\lie S_m$-modules, we see that as a $\lie S_m$-module, 
$(M_{\loc})_\mu$ is isomorphic to
 \begin{equation}\label{decomm} \left( \sum_\tau s_\mu(\tau, V) S(\tau)\right)\otimes \left(\sum_\sigma f_\sigma(u) S(\sigma)\right) = \sum_\tau \sum_\sigma s_\mu(\tau, V) f_\sigma(u) S(\tau) \otimes S(\sigma) \end{equation} 
\[ = \sum_\tau \sum_\sigma \sum_\gamma s_\mu(\tau, V) f_\sigma(u) c_{\tau, \sigma}^\gamma S(\gamma).\]
Another way to state this is that $H((M_{\loc})_\mu) = \sum_\tau \sum_\sigma \sum_\gamma s_\mu(\tau, V) f_\sigma(u) c_{\tau, \sigma}^\gamma S(\gamma).$

The Hilbert series of $B_{\loc}(\gamma, V)_\mu$ is obtained by collection the coefficient of $S(\gamma)$ in the decomposition of $(M_{\loc})_\mu$.  We know that $B_{\loc}(\gamma, V)$ is a finite-dimensional module for $\lie g$, and so is a direct sum of it's weight spaces.  Therefore the graded character of $B_{\loc}(\gamma, V)$ is given by summing over the Hilbert series for its weight spaces:
\[ \chg B_{\loc}(\gamma, V) = \sum_{\mu \in P^+}( \sum_{\tau} \sum_\sigma s_\mu(\tau, V) c_{\tau, \sigma}^\gamma f_\sigma(u) ) e(\cal O(\mu)).\]

To prove the second equality in Theorem \ref{main}, first note that by using the equality proved above, and Lemma \ref{duallem}, we have
\[ \chg B_{\loc}(\gamma^\vee, V^*)^* = \sum _{\mu \in P^+} ( \sum_\tau \sum_\sigma s_{\mu}(\tau, V^*) c_{\tau, \sigma}^{\gamma^\vee} f_\sigma(u^{-1})) e(\cal O(\mu^\vee)).\]  
It is enough to show that, for a fixed $\mu\in P^+$ and $\tau \in \part(m)$
\[ \sum_{\sigma} s_\mu(\tau, V) c_{\tau, \sigma}^\gamma f_\sigma(u) = \sum_{\sigma} s_{\mu^\vee}(\tau, V^*) c_{\tau, \sigma}^{\gamma^\vee} f_{\sigma}(u^{-1}). \]
Starting on the left, we apply Proposition \ref{thef} and Lemma \ref{thec}
\[ u^{\binom{m}{2}} \sum_{\sigma} s_\mu(\tau, V) c_{\tau, \sigma^\vee}^{\gamma^\vee} f_{\sigma^\vee}(u^{-1}).\]
Using Proposition \ref{wtspc} we get
\[ u^{\binom{m}{2}} \sum_{\sigma} s_{\mu^\vee}(\tau, V^*) c_{\tau, \sigma^\vee}^{\gamma^\vee} f_{\sigma^\vee}(u^{-1})\]
and reindexing gives the result.

\section{The case $V=V(\omega_1)$}

\subsection{Weights and Characters}

A special case of our construction is when we take $\lie g = \lie{sl}_{n+1}$ and $V$ to be the natural representation $V(\omega_1)$. In this case we can identify those dominant weights, $\mu \in P^+$, such that $(V^{\otimes m})_\mu \ne 0$ and we can explicitly describe the $\lie S_m$-structure of these weight spaces. 

The following results are well known.
\begin{prop}\label{natrep} The module corresponding to the natural representation, $V(\omega_1)$, satisfies the following properties.
\begin{enumerate}
\item It has a basis $\{ v_0,\ldots, v_n \}$, defined by $v_i \bydef x_i^- x_{i-1}^- \cdots x_1^- v_0$ for $i\ge 1$, where the $x_i$ are Chevalley basis elements.
\item Each basis vector, $v_i$, has weight, $\wt(v_i)= -\omega_i + \omega_{i+1}$, where we use the convention that $\omega_0 = \omega_{n+1}=0$.
\end{enumerate}
\end{prop}

Let $N = V^{\otimes m}$.   Let $\underline{a} = (a_0, \ldots, a_n)\in \bz_{\ge 0}^{n+1}$ be such that $\sum a_i = m$ and define $v^{\underline a} \bydef \otimes_{i=o}^n v_i^{a_i}$;  this element of $N$ has weight $\wt(v^{\underline a}) = \sum_{i=1}^n (a_{i-1} -a_i) \omega_i$.  It follows that $\wt v^{\underline a} \in P^+$ if and only if $a_{i+1} \le a_i $ for all $i$;  that is to say, precisely when $\underline{a}$ is a partition.   
Given such a partition $\underline a\in \part(m)$, define an $n$-tuple, $\underline a' \in \bz_{\ge 0}^n$, by $a_i'= \sum_{j=i}^n a_j$. We can now write the weight of $v^{\underline a}$ as
\[
\wt v^{\underline a} = m\omega_1 - \sum a_i' \alpha_i.
\] The following properties follow immediately from the definition.
\begin{lem}\label{lemma:weights}The weight spaces of $N$ can be described as follows
\begin{enumerate}[(i)]
\item Suppose $N_\mu\ne 0$, then there exists some $\underline a$ such that $\mu = m\omega_1 - \sum a_i' \alpha_i$, where $\underline a'$ is associated to $\underline a$ as above.  Furthermore, $\mu \in P^+$ if and only if $a_i \ge a_{i+1}$ for all $i$.
\item The space $N_\mu$ is spanned by permutations of the tensorands of the vector $v^{\underline a}$. 
\end{enumerate}
\end{lem}

\begin{prop} Fix $\underline a\in \part(m)$ as above and let $\mu$ be an arbitrary weight. The decomposition of the module $N_\mu$ is described by Kostka numbers; that is,
\[
N_\mu \cong \bigoplus_{\tau \ge \underline a} K_{\tau, \underline{a}} S(\tau).
\]
\begin{pf}
Define $Y(\underline{a})=\lie S_{a_0} \times \cdots \times \lie S_{a_n} \subset \lie S_m$, the Young subgroup associated to $\underline{a}$.  Clearly the action of $Y(\underline{a})$ fixes $v^{\underline{a}}$. Lemma \ref{lemma:weights} shows that $v^{\underline{a}}$ generates $N_\mu$ as a left $\CC\lie S_n$-module, that is
\[
N_\mu \cong \CC\lie S_n \otimes_{Y(\underline a)}  v^{\underline{a}}  = \mathrm{Ind}^{\lie S_n}_{\lie S_k} \Triv.
\]
The proposition follows by applying Young's Rule (see, for example, \cite[Corollary 4.39]{FultonHarris}) in the special case that the representation is trivial.
\end{pf}
\end{prop}

\begin{coro}
The character of $B_{\loc}(\gamma, V(\omega_1))$ can be described explicitly: 
\[
 \chg B_{\loc}(\gamma) = \sum_{\sigma,\tau,\underline a\in \part(m)\, |\, \tau \ge \underline a}  K_{\underline a, \tau} f_\sigma(u) c_{\tau, \sigma}^\gamma  e(\cal O(\mu_{\underline a})).
\]   
\end{coro}


\begin{thebibliography}{99}
\bibitem{BC} M.~Bennett and V.~Chari, \emph{Tilting modules for the current algebra of a simple Lie algebra}, Proceedings of Symposia in Pure Mathematics (86): Recent Developments in Lie Algebras, Groups and Representation Theory  2012, 75-97.

\bibitem{BC2} M.~Bennett and V.~Chari \emph{Character Formulae and Realization of Tilting Modules for $\lie{sl}_2[t]$}, arXiv:1409.4464. 

\bibitem{BCGM} M.~Bennett, V.~Chari, J.~Greenstein and N.~Manning, \emph{On homomorphisms between global Weyl modules}, Represent. Theory 2011, 15, 733-752.

\bibitem{C} V.~Chari, \emph{On the fermionic formula and the Kirillov-Reshetikhin conjecture}, Int. Math. Res. Not. 12 (2001), 629-654.

\bibitem{Chevalley} C.~Chevalley, \emph{Invariants of finite groups generated by reflections}, Amer. J. Math 77 (1955), 778-782.

\bibitem{CI}V.~Chari and B.~Ion \emph{BGG reciprocity for the current algebra}, Compositio Mathematica 151 (2015), 1265-1287.

\bibitem{CFK} V.~Chari, G.~Fourier and T.~Khandai, \emph{A categorical approach to Weyl modules}, Transform. Groups 15 (2010), no. 3, 517-549.

\bibitem{CPWeyl} V.~Chari and A.~Pressley \emph{Weyl modules for classical and quantum affine algebras}, Represent. Theory, 5 (2001), 191-223.

\bibitem{FultonHarris} Fulton, William and Harris, Joe, \emph{Representation theory}, 1991, Graduate Texts in Mathematics, Springer-Verlag, New York.

\bibitem{Humphreys} Humphreys, J. E., \emph{Reflection Groups and Coxeter Groups}, 1990, Cambridge University Press.

\bibitem{FL} B.~Feigin and S.~Loktev \emph{Multi-dimensional Weyl Modules and Symmetric Functions}, Comm. Math. Phys. 251 (2004), no. 3, 427-445.

\bibitem{Kh} A.~Kkoroshkin \emph{Highest weight categories and Macdonald Polynomials}, arXiv:1312.7053 

\bibitem{Kl} A.~Kleshchev \emph{Affine highest weight categories and affine quasihereditary algebras}, Proceedings London Mathematical Society 100 (2014), no. 4,  841-882.

\bibitem{M} A.~Moura \emph{Restricted limits of minimal affinizations}, Pacific J. M 244 (2010), 359-397
\end{thebibliography}
\end{document}